\documentclass[10pt,a4paper,twoside]{article}
\usepackage{amssymb}
\usepackage{latexsym}
\usepackage{fullpage}
\usepackage{a4,psfig,latexsym,rotating,amsmath,epsfig}
\usepackage{here}
\usepackage{amsmath}
\usepackage{amscd}
\usepackage{oldgerm}
\usepackage{fancyheadings}
\usepackage{here}
\def\In{\mathbb{I}}

\def\R{\mathbb{R}}

\def\D{\mathbb{D}}
\def\1{\mathbf{1}}
\def\:{\lrcorner}
\def\#{\sharp}

\def\a{\alpha}
\def\b{\beta}
\def\g{\gamma}
\def\d{\delta}
\def\e{\epsilon}
\def\p{\partial}

\def\o{\circ}

\def\x{\otimes}
\def\<#1,#2>{\langle#1,\,#2\rangle}

\def\G{\Gamma}

\def\S{\mathbb{S}\,}

\def\qed{\ensuremath{\quad\Box\quad}}

\def\inv#1{\raise.1em\hbox to 0pt{$^{-1}$\hss}_{#1}\;}

\newtheorem{Theorem}{Theorem}
\newtheorem{Lemma}[Theorem]{Lemma}

\newtheorem{Definition}[Theorem]{Definition}

\title{The Cauchy Problem of Lorentzian Minimal Surfaces in Globally Hyperbolic  Manifolds}

\begin{document}

\maketitle

\begin{center}

Olaf M\"uller\\Max-Planck-Institute for mathematics in the sciences\\Inselstraße 22-26\\D-04105 Leipzig, Europe

\end{center}

\begin{abstract}
\noindent In this note a proof is given for global existence and uniqueness of minimal Lorentzian surface maps from a cylinder into globally hyperbolic Lorentzian manifolds for given initial values up to the first derivatives.

\end{abstract}

\section{Introduction}

\noindent In the search for a unified field theory of all interactions it is widely believed that a perturbative expansion of such a theory could be provided by String Theory. Its classical action, at least in the case of Closed Bosonic String Theory, is the area of a surface mapped into spacetime. As the latter one is supposed to have Lorentz signature the idea of a one-dimensional object moving through the space suggests that one should consider only mappings that induce a Lorentzian metric on the surface to get a consistent definition of area. While harmonic mappings of Lorentzian surfaces into Riemannian manifolds are a well-examined object since the pioneering work of Gu (\cite{gu}), there seems to be no comparable global result in the double-Lorentzian case in the literature. All gradient estimates (a main tool in Gu's paper) fail in this case because of the presence of the null cone, so one has to apply different methods. The result given by Shatah and Struwe in ~\cite{ss} is a local one without time estimates. The main result of this section is that the minimal surface problem is globally well-posed if we assume global hyperbolicity for the target manifold (which replaces in some sense the tacit assumption of bounded geometry in Gu's paper). The solution will exist {\em locally} on the string worldsheet, but {\em globally} in the target space in the sense that the image will intersect {\em any} time slice in a standard way. The results extend earlier results of Thomas Deck (\cite{de}) for static spacetimes.

\bigskip

\noindent The minimal surface problem for a {\em Riemannian} target space was solved by J.Douglas in the case of a flat target space (\cite{do}) and by C.B.Morrey in the case of a general Riemannian manifold as target space (\cite{mo}). The problem for flat Lorentzian target space is trivial as it consists basically in solving the linear wave equation. But in the case of non-flat Lorentzian target manifolds left- and rightmovers in the language of physics do not decouple any more such that the problem is more involved. We are going now to adress this remaining gap in the case of a globally hyperbolic target manifold giving an affirmative answer.

\section{General facts}

\noindent  Let $(M,g)$ be a Lorentzian manifold of dimension $n$. We first recall some well-known facts about minimal surfaces. 

\begin{Definition}
The area of a map $y : \S^1 \times \In \rightarrow M$ for which the pulled-back metric is Lorentzian everywhere is defined to be
$$A(y) := \int_{(\S^1 \times \In, y^* g)} 1 .$$
\noindent If $y$ is a critical point of the area functional in the space of smooth maps it is called minimal surface.
\end{Definition}

\noindent This functional has a lot of symmetries: it is invariant under diffeomorphisms of $N$.

\begin{Definition}
Let $(N,h)$ be a two-dimensional Lorentzian manifold. A map $y : S \rightarrow M$ is called {\em wave map} iff it satisfies the following differential equation:

$$ tr \nabla d y = 0 $$

\noindent where the trace is understood with respect to the metric $h$ on $N$, that means in coordinates $t,x$ which are orthogonal at a point $p \in N$:

$$ (\partial_t)^2 y ^{\a}- (\partial_x)^2 y ^ {\a} + \G ^{\a}_{\b\g} (y) ( \partial_t y^{\b} \partial_t y^{\g} -  \partial_x y^{\b} \partial_x y^{\g}) = 0$$ 

at $p$.
\end{Definition}

\noindent Equivalently, we can write the equations in characteristic coordinates  $\xi := x+t , \eta := x-t$ on $\Delta$ as 

$$\partial_{\xi} \partial_{\eta} y^{\a} +  \G ^{\a}_{\b\g} (y)  \partial_{\xi} y^{\b} \partial_{\eta} y^{\g} =0 .$$

\begin{Theorem}
$y$ is a minimal surface iff it is a wave map with respect to $y^* g$. 
\end{Theorem}

\begin{Definition}
For a map $y: (N,h) \rightarrow (M,g)$ the energy $E(y)$ is defined as

$$E(y) := \int_{(N,h)} \langle df, df \rangle  $$

where the norm is the operator norm. If $y$ is a critical point of $E$ in the space of smooth maps it is called minimal-energy surface.

\end{Definition}

The energy is invariant under conformal changes of $h$. The next theorem is an immediate consequence from the Euler Lagrange equations:

\begin{Theorem}
$y$ is a minimal energy surface iff it is a wave map and $y^* g$ is conformally equivalent to $H$.
\end{Theorem}

Thus, minimal energy surfaces are minimal surfaces, and while area and energy differ as functionals they have the same critical points. By diffeomorphism invariance and by the result of Kulkarni (\cite{ku}) that every Lorentzian metric on a cylinder wit Lorentzian ends can be linked by a diffeomorphism to a conformal multiple of the standard metric, we get

\begin{Theorem}
$y: \S^1 \times [ 0 , 1 ] $ is a minimal surface if it is a wave map w.r.t. the standard metric and if $y^*g$ is conformally equivalent to the standard metric.

\end{Theorem}

We will first focus on the first property: that $y$ be a wave map.

\section{Wave maps. Local existence and uniqueness}

\noindent We will look for wave maps $y$ from some subsets of the upper half Minkowski plane to $M$ with some initial values on the x-axis. Throughout the proof we will work in characteristic coordinates $\xi := \frac{t+x}{2}$, $\eta := \frac{t-x}{2}$ on $\R^{1,1}$. Let $\frac{\partial y ^{\a} }{ \partial \xi} = u^{\a} ,  \frac{\partial y ^{\a} }{ \partial \eta} = v^{\a}$. 
We will need an auxiliary Riemannian metric $h = \langle \cdot , \cdot \rangle^+$ and an associated norm $\vert \vert . \vert \vert$ on $M$. Define, for every $p \in M$, $R(p)$ as the injectivity radius of $h$ at $p$, fix coordinate patches at every $p$ associated to these balls, let $\Gamma^{\a}_{\b \gamma} $ be the Christoffel symbols of the Levi-Civita connection of $g$ in the chosen coordinates, let $G(p)$ be the maximum of the operator norm (w.r.t. h) of $\Gamma_{\a \b}^{\gamma} $ in the coordinate patch around $p$, similarly $\underline{u} , \underline{v} := \vert \vert u  \vert \vert_{C^1( \{ t=0 \} )} , \vert \vert v  \vert \vert_{C^1( \{ t=0 \} )}$, respectively.

\begin{Theorem}
\label{wave}
Let $M$ be a Lorentzian manifold, let $k=(k_0,k_1) : \R \rightarrow TM$ be a smooth curve of bounded derivatives, i.e. $ |\partial_s k_i \vert \leq c_{|s|} \in \R$ for all $s$ and $\langle \partial_x k_0 , \partial_x k_0 \rangle  \geq 0 \geq  \langle k_1 , k_1 \rangle $ everywhere. Let $G(p)$ be bounded by $G(\d)$ in a neighborhood $B_{\d} (k_0)$ of $im(k_0)$. Then there is a unique smooth map from a neighbourhood of the x-axis to $M$  which is a wave map and whose restriction and normal derivative at $\{ t=0 \}$ correspond to the given curve, i.e. $y(0,x) = k_0 (x)$, $\partial_t y(0,x) = k_1(x)$.
Let $R_k:= \inf_{p \in im(k_0)} R(p)$, $L_k := \min \{ \frac{R_k}{5} , \frac{1}{11 G(\d)} , \frac{\d}{5} \}$, then the neighbourhood can be chosen as a small strip $y: \{ t \in [ 0 , \frac{1}{\sqrt{2}} l ] \} \rightarrow M $ with 
$ l :=  \frac{L_k}{  (\underline{u} + \underline{v} ) }$.
Note that if $k$ is a closed curve in $TM$ then $l$ is bounded from zero.

\end{Theorem}

\noindent The theorem will be proved as an easy corollary of the same statement for small characteristic triangles based on the boundary curve $ \{ t=0 \}$:

\begin{Lemma}
Let $M$ be a Lorentzian manifold, let $k=(k_0,k_1) : \R \rightarrow TM$ be as above, $l$ as above. Then for each characteristic triangle $\Delta$ of base length $\leq l$ there is a unique smooth map $y: \Delta \rightarrow M$ which is a wave map and whose restriction and normal derivative at $\{ t=0 \}$ correspond to the given curve, i.e. $y(0,x) = k_0 (x)$, $\partial_t y(0,x) = k_1(x)$.
\end{Lemma}

\bigskip

\noindent {\em Proof.} First note that the choice of $l$ as above assures that the image of each interval of length $l$ is contained in some coordinate patch as $l \leq   \frac{R_k}{2|| \partial_x k_0 ||^+} $. Let $\Delta$ be the triangular region spanned by the interval of length $\leq l$ (which we call the base side $a$ of the triangle) and the characteristic lines beginning at the endpoints of the interval. For every point $p$ in $\Delta$ we define four special lines: The lines $c_{p, \eta}$ resp. $c_{p, \xi}$ are just the characteristic lines along the constant vector fields $\partial_{\eta}$ resp. $\partial_{\xi}$ ending at $p$ and beginning at the basis of the triangle at points we call $p'$ resp. $p''$ while the lines $\underline{c}_{p, \eta}$ resp. $\underline{c}_{p, \xi}$ are the lines on the base side from the left endpoint of the base side to $p'$ resp. $p''$. We use a way of splitting the differential equation similar to the one in Gu's paper (\cite{gu}) and consider the following system of first-order ordinary differential equations of maps $y,z,u,v, \hat{u} , \hat{v} $ from $\Delta$ to $\R ^n$ :

\begin{align}
\label{vier}
 \frac{ \partial u ^ {\a} }{\partial \eta } + \G^{\a}_{\b \g} (z) u^{\b} \hat{v}^{\g} = 0 , \qquad  &  \frac{ \partial \hat{u} ^ {\a} }{\partial \eta } + \G^{\a}_{\b \g} (y) \hat{u}^{\b} v^{\g} = 0, \nonumber\\ 
 \frac{ \partial v ^ {\a} }{\partial \xi } + \G^{\a}_{\b \g} (z) \hat{v}^{\b} u^{\g} = 0 , \qquad  & \frac{ \partial \hat{v} ^ {\a} }{\partial \xi } + \G^{\a}_{\b \g} (y) v^{\b} \hat{u}^{\g} = 0, \\
 \frac{\partial y ^{\a} }{ \partial \xi} = u^{\a} ,  \frac{\partial y ^{\a} }{ \partial \eta} = v^{\a}, \qquad & \frac{\partial z ^{\a} }{ \partial \xi} = \hat{u}^{\a} ,  \frac{\partial z ^{\a} }{ \partial \eta} = \hat{v}^{\a} \nonumber
\end{align}

\noindent with initial conditions 

$$y^{\a} (0,x) = z^{\a} (0,x) = k ^{\a}_{0} (x),$$

$$u ^{\a} (0,x) = \hat{u} ^{\a} (0,x)  = \frac{\partial k^{\a}_0 (x) }{\partial x} + k^{\a}_1 (x) ,\ \  v ^{\a} (0,x) = \hat{v} ^{\a} (0,x) = - \frac{\partial k^{\a}_0 (x) }{\partial x} + k^{\a}_1 (x).$$

\noindent Note that for example the first equation of the system (\ref{vier}) corresponds to the invariant equation $\nabla_{\eta} u = 0$ (with $\nabla$ denoting the covariant derivative), an equation of parallel transport along $\eta$-lines. A solution of this system will give $y=z$ (because of uniqueness and symmetry of the equations under $y \leftrightarrow z, u \leftrightarrow \hat{u}, v \leftrightarrow \hat{v}$) and therefore be a solution of the original problem. At this point we can forget about the Lorentzian metric $g$ that does not appear at all in the equations and use the Riemannian metric $h$ as long as we still keep the Christoffel symbols of $g$.

\bigskip

\noindent Now consider the following iteration procedure whose fix points are exactly the solutions of ( \ref{vier}): Let $(y_0,z_0,u_0,v_0, \hat{u}_0, \hat{v} _0)$ be any system of maps satisfying the respective initial conditions and 

$$ \frac{\partial y_0^{\a}}{\partial \xi} = u_0^{\a} ,\  \frac{\partial y_0 ^{\a} }{ \partial \eta} = v_0^{\a}, \ \frac{\partial z_0 ^{\a} }{ \partial \xi} = \hat{u}_0^{\a}, \ \frac{\partial z_0^{\a}}{\partial \eta} = \hat{v}_0^{\a} $$

\noindent Now the iteration process 

$(y_m,z_m,u_m,v_m, \hat{u}_m, \hat{v}_m) \mapsto (y_{m+1},z_{m+1},u_{m+1},v_{m+1}, \hat{u}_{m+1}, \hat{v}_{m+1}) $ is defined by

\begin{align}
\label{iter}
\frac{\partial u_{m+1}^{\a} }{\partial \eta} + \G^{\a}_{\b \g} (z_m) u_{m+1}^{\b} \hat{v}_m^{\g} = 0 , \qquad & \frac{\partial \hat{u}_{m+1}^{\a} }{\partial \eta} + \G^{\a}_{\b \g} (y_m) \hat{u}_{m+1}^{\b} v_m^{\g} = 0 \nonumber \\
 \frac{\partial v_{m+1}^{\a} }{\partial \xi} + \G^{\a}_{\b \g} (z_m) \hat{v}_{m+1}^{\b} u_m^{\g} = 0, \qquad & \frac{\partial \hat{v}_{m+1}^{\a} }{\partial \xi} + \G^{\a}_{\b \g} (y_m) v_{m+1}^{\b} \hat{u}_m^{\g} = 0 \\
 \frac{\partial y_{m+1}^{\a} }{\partial \xi} = u_{m+1}^{\a} , \ &  \frac{\partial z_{m+1}^{\a} }{\partial \eta} = \hat{v}_{m+1}^{\a} \nonumber
\end{align}

\noindent Note that a fix point of this iteration procedure will also be $(y \leftrightarrow z,  u \leftrightarrow \hat{u}, v \leftrightarrow \hat{v})$-symmetric because of the initial conditions and 

$$\frac{\partial v}{\partial \xi} = - \Gamma (z) \hat{v} u = \frac{\partial u}{\partial \eta} = \frac{\partial y}{\partial \eta \partial \xi}$$

\noindent which yields $ v = \frac{\partial y}{\partial \eta}$ and the other missing equation of system (\ref{vier}). It will therefore be a solution of the original problem. Now we consider only  maps $\Delta \rightarrow M$ which have the origin of our local coordinate system as the image of the left endpoint of the base side, thus we can forget about $y,z$ just by inserting

$$y_{m} (p) =   \int_0^{l(c_{p, \xi})} u_m (c_{p, \xi} (t)) dt
 + \int_0^{l(\underline{c}_{p,\xi})} (\hat{u}_m - \hat{v}_m) (\underline{c}_{p,\xi} (t) ) $$

$$z_{m} (p) =   \int_0^{l(c_{p, \eta })} \hat{v}_m (c_{p, \eta} (t)) dt
 + \int_0^{l(\underline{c}_{p,\eta})} (\hat{u}_m - \hat{v}_m) (\underline{c}_{p,\eta} (t) ) $$

\noindent (every curve, defined as above, is parametrized by arc length in Euclidean $\R^2$, $l$ denotes the respective lengths). As stated above, this iteration procedure is well-defined and meaningful in the space of maps with different images of the base sides but one and the same image of the left endpoint of the base side. Nevertheless the initial data remain fixed during the iteration procedure, $u_{m+1} \vert_{\{t=0 \} } =u_{m} \vert_{\{t=0 \} }$ etc.; the corresponding operator taking $r_m := (u_m,v_m, \hat{u}_m, \hat{v}_m)$ to $r_{m+1} := (u_{m+1},v_{m+1}, \hat{u}_{m+1}, \hat{v}_{m+1})$ we denote by $\Phi  = (\Phi_1 , \Phi_2 , \Phi_3 , \Phi_4 )$. It is mapping $C^1$ maps to $C^1$ maps which is seen trivially from the defining equation for half of the coordinates, for the other half of coordinates use the usual theorems about differentiable dependence of ODE systems on parameters as in \cite{wa} , II. 13. V. From now on, we will omit the sequential index of the $r_n$; a subindex will now be related to $r = (r_1, r_2, r_3, r_4)$.

\bigskip

\noindent Now we have a special solution of the system (\ref{vier}) to compare with, namely the one with $ \partial _t \vert_{ t=0} y = \pm  \partial _x \vert_{ t=0} y$, i.e. with $u=0=\hat{u}$ or $v=0=\hat{v}$ everywhere, we decide for $u=0=\hat{u}$. It is easy to see that this is a solution (for a different initial value curve $\tilde{k}$, but still $ \tilde{k}_0 = k_0 $), we call it $s:=(0,V,0,V)$ (with $V:=(v-u) | _{ \{ t=0 \} }$). Its image is one-dimensional, namely the chosen interval of the curve $k_0$; as it is a solution, it is a fix point of $\Phi$:  $\Phi ( (0,V,0,V)) = (0,V,0,V)$. Moreover, as we assumed $ l \leq  \frac{R_k}{ 2 ||u-v||_{\{t=0\}}^+ }$, $l$ is so small that the whole image of the base of the triangle is contained in some coordinate patch with distance $\geq \frac{R_k}{2}$ from its boundary, thus we have a ball-shaped coordinate patch $B$ containing the image of $s$ with

\begin{equation}
\label{hehe}
dist(im (s) , \partial B ) \leq \frac{R_k}{2}
\end{equation}

\noindent Now consider the affine subspace $M (r) := \{ f \in C^1 ( \Delta ,( \R ^n )^4 ) | \ f |_a = r |_a \} $, define seminorms  $||.||_j$ of maps from $\Delta$ to $(\R^n)^4$ by

$$|| f ||_0 := max_{i=1,...4} \{  || f_i || _{C^0 (\Delta)}  \} $$

$$|| f ||_j := max_{i=1,...4} \{  || d^j f_i || _{C^0 (\Delta)}  \} $$

\noindent Note that in the case $j=1$ the corresponding distance to the special solution (or to any other map in $M(r)$) is then a metric on $M (r)$ as $||f-g||_1 = 0$ means that $f$ and $g$ differ by an additive constant.

\bigskip

\noindent Now, for given $ K \in \R$, $r,s \in C^{\infty}(\Delta, \R^M)$ let us define 

$$M_K (r,s) := \{f \in M (r) | \ ||f-s||_1 \leq K \} .$$

\noindent Note that in general $s$ does not need to be in $M(r)$. For technical reasons, we consider the reflection $s$ at the symmetry axis of the characteristic triangle and the map $s' :=(\1, s , \1, s)$. We define $\hat{\Phi} : = s' \o \Phi \o s'$. 

\begin{Lemma}
Recall $\underline{u} := ||u||^{C^1}_{\{t=0\}} , \underline{v} := ||v||^{C^1}_{\{t=0\}}$. Define $ K := \frac{3}{l} (\underline{u} + \underline{v}) $, then

(i) If $r \in C^1$, then $\partial_{\eta} \partial_{\eta} (\hat{\Phi} (r)) $ and $\partial_{\eta} \partial_{\xi} (\hat{\Phi} (r))$ exist and are continuous.

(ii) $\hat{\Phi} (r) |_a = r|_a$, i.e. the values of the maps on the base side remain fixed.

(iii) Let $s$ be the special solution to $r$ as above, then $||r-s||_1 \leq K $ implies that $|| \hat{\Phi} (r) -s ||_1 \leq K$.

\end{Lemma}

\noindent {\em Proof of the lemma.} (i) and (ii) follow easily from the defining equations. For (iii) use that on the base side $a$, $s$ and $r$ differ at most by the length $(\underline{u} + \underline{v})$, so $||r-s||_1 \leq K $ implies $||r-s||_0 \leq Kl + (\underline{u} + \underline{v}) = 4  (\underline{u} + \underline{v}) =: K'$.

As $||r-s||_0 < K'$  we have the following estimates:

$$ || \Gamma^{\a}_{\b \gamma} (y_m) \cdot v^{\gamma} || \leq G(\d) \cdot (\underline{u} + \underline{v} + ||r-s||_0 ) \leq 5G(\d) (\underline{u} + \underline{v})$$

\noindent and thus standard ODE estimates (~\cite{wa} , 2.12.V.) give

$$ || \hat{\Phi}_1 (r) ||_{C^0}, || \hat{\Phi}_2 (r) -V ||_{C^0}, ||\hat{\Phi}_3 (r) ||_{C^0} , ||\hat{ \Phi} _4 (r) -V ||_{C^0} \leq (\underline{u} + \underline{v}) e^{5G l(\underline{u} + \underline{v})} \leq K' $$

\bigskip

\noindent Now we want to get the same result for $K$ instead of $K'$ and on the $C^1$ level. For half of the coordinates this follows easily from the defining equations and the $C^0$ estimates, as e.g. 

$$ ||\frac{\partial u_{m+1}}{\partial \eta} ||_{C^0} = ||- \G (z_m) \cdot u_{m+1} \hat{v}_m||_{C^0} \leq  (K')^2 G \leq \frac{3}{8} \frac{K'}{l} = \frac{K}{2}$$

\noindent as $l \leq \frac{3}{32 (\underline{u} + \underline{v}) G}$. The other half we get by applying again the quoted standard ODE theorem (\cite{wa}, 2.13.V. and 2.12.V.) to $\partial_x r$ ($r$ considered to be fixed on the corresponding characteristic line)  where the norm of the Lipschitz coefficient is again $\leq G \cdot K'$. Then the comparison with the fix point $s$ of $\hat{\Phi}$ yields $||\partial_x (\hat{\Phi} (r)) - \partial_x s ||_0 \leq \g e^{L l}$ where $\g$ is the distance $||\partial_x r- \partial_x s||_{C^0 (a)}$ of the maps restricted on $a$ which is less or equal to $\frac{K}{4}$ while $L$ is the Lipschitz coefficient which can be taken to be $G \cdot K'$, so 

$$ ||\partial_x r- \partial_x s||_0 \leq \frac{K}{4} e ^{GK'l} < \frac{K}{2}$$

\noindent which completes the proof of the lemma   \qed

\noindent In short, the lemma above tells us that $\Phi (M_K (r,s)) \subset M_K (r,s)$ (and that differentiasbility improves a little bit) if we choose $K$ as above.

\begin{Lemma}
The set $A:=\Phi (conv (\Phi( M_K (r,s))))$ has compact closure in $C^1 (\Delta, \R^N)$.    

\end{Lemma}

{\em Proof.} According to Arzela-Ascoli, a set in $C^1$ is precompact if and only if it is $C^1$-bounded and all of its first derivatives are equicontinuous (Lipschitz-bounded). But following the lemma above, $\Phi (M_K (r,s)) \subset M_K (r,s)$, $conv(\Phi (M_K (r,s))) \subset M_K (r,s)$ as the norms are convex, and $\Phi (conv(\Phi (M_K (r,s)))) \subset M_K (r,s)$ following the lemma again. Therefore $A$ is bounded in $C^1$.

\noindent Now $u, v, \hat{u}, \hat{v}$ are equicontinuous as they are $C^1$-bounded. Analogously for $\partial_{\eta} u$, $\partial_{\xi} v$, $ \partial_{\eta} \hat{u}$, $\partial_{\xi} \hat{v}$ (look at the defining equations). Thus it remains to be shown that $\partial_{\xi} u$, $\partial_{\eta} v$, $ \partial_{\xi} \hat{u}$, $\partial_{\eta} \hat{v}$ are equicontinuous in $A$. First one can establish equicontinuity in one half of the directions, i.e. show that $\partial _{\xi} u (p) - \partial _{\xi} u (p + \tau \cdot \eta) \leq C \cdot \tau$. This can be proven by the usual ODE theorem about differential dependence on the parameter $\eta$ and then performing the admissible formal differentiation 

$$\partial_{\eta} (\partial_{\xi} u_2) = \partial_{\xi} (\partial_{\eta} u_2) = \partial_{\xi} ( - \Gamma(y_1) u_2 v_1)$$

\noindent which shows equicontinuity w.r.t. one half of the variables.

\noindent  The other half one can get by differentiating the defining equation as above

\begin{equation}
\label{krass}
\frac{\partial}{\partial \eta} \frac{\partial u_2^{\a}}{\partial \xi} = - ( \frac{\p}{\p \xi} \Gamma^{\a}_{\b \g} (z_1) ) \hat{v}_1^{\g} u_2^{\b}  -    \Gamma^{\a}_{\b \g} (z_1) )(\frac{\p}{\p \xi} \hat{v}_1^{\g}) u_2^{\b}   -   \Gamma^{\a}_{\b \g} (z_1)  \hat{v}_1^{\g} ( \frac{\p}{\p \xi} u_2^{\b}) =: f ( \eta , \frac{\p u_2}{\p \xi }) ,
\end{equation} 

\noindent understood as the defining ODE for $\frac{\p u}{\p \xi}$ along a fixed $\eta$-line for given $u$. Now we compare the $\eta$-line with a second one of distance $\e$ for which all quantities are denoted by $u^{(\e)}, \hat{v}^{(\e)}, z^{(\e)}$:

$$f^{(\e)} ( \eta , \frac{\p u_2}{\p \xi } ):=  - ( \frac{\p}{\p \xi} \Gamma^{\a}_{\b \g} (z_1^{(\e)}) ) \hat{v}_1^{(\e) \g} u_2^{ (\e) \b}  -    \Gamma^{\a}_{\b \g} (z_1^{(\e)}) (\frac{\p}{\p \xi} \hat{v}_1^{ (\e) \g}) u_2^{(\e) \b}   -   \Gamma^{\a}_{\b \g} (z_1^{(\e)})  \hat{v}_1^{(\e) \g} ( \frac{\p}{\p \xi} u_2^{(\e) \b})$$

\noindent The equicontinuity will be shown by the theorem about continuous dependence of ODEs on a real parameter (\cite{wa}, II.12.VI). To do this, we have to prove that for any $\d >0$ there is an $\e >0$ such that

$$f^{(\e)} (\eta, \frac{\p u_2}{\p \xi} ) -f (\eta, \frac{\p u_2}{\p \xi} ) \leq \d .$$
 
\noindent So we have to consider 

\begin{align}
- \frac{\p}{\p \xi} \Gamma^{\a}_{\b \g} (z_1) ) \hat{v}_1^{ \g} u_2^{  \b} + \frac{\p}{\p \xi} \Gamma^{\a}_{\b \g} (z_1^{(\e)}) ) \hat{v}_1^{(\e) \g} u_2^{ (\e) \b} \nonumber\\
- \Gamma^{\a}_{\b \g} (z_1) )(\frac{\p}{\p \xi} \hat{v}_1^{ \g}) u_2^{ \b} + \Gamma^{\a}_{\b \g} (z_1^{(\e)}) )(\frac{\p}{\p \xi} \hat{v}_1^{ (\e) \g}) u_2^{(\e) \b} \nonumber\\
- \Gamma^{\a}_{\b \g} (z_1)  \hat{v}_1^{ \g} ( \frac{\p}{\p \xi} u_2^{ \b}) + \Gamma^{\a}_{\b \g} (z_1^{(\e)})  \hat{v}_1^{(\e) \g} ( \frac{\p}{\p \xi} u_2^{(\e) \b}) \nonumber
\end{align}

\noindent The first line is quite easily shown to be less than $3GK^3 \e$ as every difference between terms with $( \e)$ and such without is less than $K \e$ or $G \e$, equally the last line. Most subtle is the second one as it seems we have not enough control over $\frac{\p}{\p \xi} \hat{v}_1^{ \g} - \frac{\p}{\p \xi} \hat{v}_1^{ (\e) \g}$. But on $\Phi (M_K(r,s))$ this term is bounded linearly in $\e$ as $\frac{\partial \hat{v}_1^{\a} }{\partial \xi} = - \G^{\a}_{\b \g} (y_0) v_1^{\b} \hat{u}_0^{\g}$. This Lipschitz bound survives the convex closure as the $C^1$ Lipschitz norm is continuous in $C^1$ and convex. Thus the second line is equally bounded in linear terms of $\e$ and we are done   \qed

\bigskip

\noindent Therefore we can apply the Schauder-Tychonoff fix point theorem (\cite{gd}, p.148) to the operator $\hat{\Phi} $ restricted to the closed convex set $conv (\Phi (M_K (r,s))) $ and prove the existence of a fix point in $C^1$ in $(u, v, \hat{u}, \hat{v})$, i.e. of a $C^2$-fix point in $(y,z)$. Below we will prove uniqueness of such fix points. Regularity of the solution is easily seen by succesive aplication of the quoted ODE theorem (\cite{wa}, 2.13.V. and 2.12.V.) to the solution (the single equations of the system \ref{vier} considered as solutions to ODEs)    \qed

\bigskip

\begin{Lemma}
$C^2$-solutions $(y,z)$ of (\ref{vier}) are uniquely determined by their initial values on the base $a$ of the triangle. 
\end{Lemma}

\noindent {\em Proof.} The proof will use energy estimates for symmetric hyperbolic systems, thus we want to bring the system into this form. To this aim we rewrite the equations in a more symmetric way, with only $y,z$ as variables:

\begin{align}
\label{uni}
\frac{\partial}{\partial \eta} \frac{\partial}{\partial \xi} y^{\a} + \G^{\a}_{\b \g} (z) \frac{\partial y^{\b}}{\partial \xi} \frac{\partial z^{\g}}{\partial \eta} &=0 \\
\frac{\partial}{\partial \xi} \frac{\partial}{\partial \eta} z^{\a} + \G^{\a}_{\b \g} (y) \frac{\partial z^{\b}}{\partial \xi} \frac{\partial y^{\g}}{\partial \eta} &=0 \nonumber
\end{align}

\noindent with initial conditions $y((x,0)) = k_0 (x) = z((x,0))$ and $\partial_t y((x,0)) = k_1(x) = \partial_t z((x,0))$. Although this looks like only a part of the system (\ref{vier}), namely the first and the last equation, it is actually equivalent to it, because if $(y,z)$ is a solution to (\ref{uni}), $(z,y)$ is as well, so uniqueness of the solution implies that then $y=z$ holds which causes equivalence to the system (\ref{vier}). We arrange the system still a bit differently:

$$\frac{\partial}{\partial x} \frac{\partial}{\partial x} y^{\a}
-\frac{\partial}{\partial t} \frac{\partial}{\partial t} y^{\a}  + \G^{\a}_{\b \g} (z) \frac{\partial y^{\b}}{\partial x} \frac{\partial z^{\g}}{\partial x} 
- \G^{\a}_{\b \g} (z) \frac{\partial y^{\b}}{\partial t} \frac{\partial z^{\g}}{\partial t}
=0$$ 

$$\frac{\partial}{\partial x} \frac{\partial}{\partial x} z^{\a}
-\frac{\partial}{\partial t} \frac{\partial}{\partial t} z^{\a}  
+ \G^{\a}_{\b \g} (y) \frac{\partial z^{\b}}{\partial x} \frac{\partial y^{\g}}{\partial x} 
- \G^{\a}_{\b \g} (y) \frac{\partial z^{\b}}{\partial t} \frac{\partial y^{\g}}{\partial t}
=0$$ 

\noindent With the notation $U := (y,z) $ we can define the functionals 

$$F^{\a} (U) := \frac{\partial}{\partial x} \frac{\partial}{\partial x} y^{\a}
-\frac{\partial}{\partial t} \frac{\partial}{\partial t} y^{\a}  + \G^{\a}_{\b \g} (z) \frac{\partial y^{\b}}{\partial x} \frac{\partial z^{\g}}{\partial x} 
- \G^{\a}_{\b \g} (z) \frac{\partial y^{\b}}{\partial t} \frac{\partial z^{\g}}{\partial t} \  {\rm for} \   \a \leq n-1$$

$$F^{\a} (U) := \frac{\partial}{\partial x} \frac{\partial}{\partial x} z^{\a}
-\frac{\partial}{\partial t} \frac{\partial}{\partial t} z^{\a}  + \G^{\a}_{\b \g} (y) \frac{\partial z^{\b}}{\partial x} \frac{\partial y^{\g}}{\partial x} 
- \G^{\a}_{\b \g} (y) \frac{\partial z^{\b}}{\partial t} \frac{\partial y^{\g}}{\partial t}
\ {\rm for } \  \a \geq n$$

\noindent Obviously $U$ is a solution if and only if the $\R^{2n}$-valued functional $F$ vanishes. Now let $N=2n$, $X=(x,t)$, then the PDE system is now in the form

$F(X,U(X), DU(X), D^2 U(X)) =0$

$F: S := \Delta \times \R^N \times (\R^2 \x \R^N) \times (\R^2 \x _s \R^2) \x \R^N \rightarrow \R^N$
 
\noindent where $\x _s$ denotes the symmetrical tensor product. An element of $S$ will be written in the form $(X,z,p,r)$. Now, following ~\cite{jo} we want to construct a suitable quasilinear operator $L: S \rightarrow \R^N$ with 

$L^{\alpha} (U_1 - U_0) = F^{\a}(U_1) - F^{\a} (U_0) .$

\noindent This is done by defining

$$V^{\a} := U_1^{\a} - U_0^{\a}$$

$$U_{\tau}^{\a} := \tau U_1^{\a} + (1- \tau ) U_0^{\a}$$

$$a_{ij \b}^{\a} (X) := \int_0^1 \frac{\partial F^{\a} }{\partial r_{ij}^{\b}} (X,U_{\tau}(X), DU_{\tau}(X), D^2 U_{\tau} (X) ) d \tau$$

$$ b_{i \b}^{\a}(X) := \int_0^1 \frac{\partial F^{\a}}{\partial p_i^{\b}} (X,U_{\tau}(X), DU_{\tau}(X), D^2 U_{\tau} (X) ) d \tau$$

$$c_{\b}^{\a} (X) := \int_0^1 \frac{\partial F^{\a}}{\partial z^{\b}} (X,U_{\tau}(X), DU_{\tau}(X), D^2 U_{\tau} (X) ) d \tau$$

and

$$L V := \sum_{i,j=1}^2 \sum_{\b=1}^N a_{ij \b }^{\a} V_{x^i x^j}^{\b} + \sum_{i=1}^2 \sum_{\b=1}^N b_{i \b}^{\a} V_{x^i}^{\b} + \sum_{\b=1}^N c_{\b}^{\a} V^{\b} ,$$

\noindent everything taken at $X=(x,t)$. All coefficients are sufficiently smooth and uniformly bounded, $a_{ij \b}^{\a}= \d_{ij} \e _i \d_{\b}^{\a}$. We have to show that if $LV=0$, $V |_{ \{ t=0 \} } = 0$ then $V=0$ everywhere. Now using techniques shown in ~\cite{fj} we can define $w \in \R ^3 \x \R ^N$, $w_1 = V_x , w_2 = V_t, w_3 = V$ and transform the system for the very last time into

$$\partial_t w_1^{\a} - \partial_x w_2^{\a} = 0$$

$$\partial_t w_2^{\a} - \partial_x w_1^{\a} - \sum_{\b =0}^{N-1} (b_{1 \b}^{\a} w_1^{\b} +   b_{2 \b}^{\a} w_2^{\b}) -  \sum_{\b =0}^{N-1} c_{\b}^{\a} w_3^{\b} = 0$$

$$\partial_t w_3^{\a} - \partial_x w_3^{\a} - w_2 + w_1 =0$$

or in short

$$Lw = A \partial_t w + A^1  \partial_x w + B(x,t,w) w = 0$$

\noindent with $A= \1 _{\R^3 \x \R^N}$, 

$$ A^1 = - \left( \begin{array}{ccc} 
0 & \1_{\R^N} & 0 \\ \1_{\R^N} & 0 & 0 \\ 0 & 0 & \1_{\R^N} 
\end{array} \right)  , \qquad B = \left( \begin{array}{ccc} 
0 & 0 & 0 \\ 
\sum_{\b} b_{1 \b}^{\a} & \sum_{\b} b_{2 \b}^{\a} & - \sum_b c_{\b}^{\a}       \\ \1_{\R^N} & -\1_{\R^N} & 0 
\end{array} \right) $$

\noindent Again $B$ is sufficiently smooth and uniformly bounded. We multiply the above equation by $w^T$ from the left and get by symmetry of the $A$-matrices

$$\partial_t (w^T A w) + \partial_x (w^T A^1 w) = 2w^T B w$$
 
\noindent and integrating over a domain $R$ and applying the divergence theorem

$$\int_{\partial R} w^T (A \frac{dt}{d \nu} + A^1 \frac{dx}{d \nu} ) w = 2 \int_R w^TBw$$

\noindent where $\nu$ is the outer normal vector of the region. Now, if we take $R_s = \Delta \cap \{ t \leq s \} $, on the base side $w$ vanishes by assumption, the terms coming from the left and right boundary are easily seen to be positive definite, so if $\frac{K}{2}$ is an upper estimate for the operator norm of $B$ we get

$$E(s) := \int_{ \{ t= s \} } w^T w \leq 2 \int_{ R_s} w^TBw \leq K \int_{ R_s} w^T w ,$$

\noindent thus $E(s) \leq K \int_0^s E(t) dt$, therefore $\partial_s E(s) \leq K e(s)$ and $ \partial_s( e^{-Ks} \int_0^s E(t) dt) \leq 0 $. But $e^{-Ks} \int_0^s E(t) dt$ is positive by definition and vanishes at $s=0$, so it vanishes everywhere, thus $w$ as well. This completes the argument  \qed

\bigskip

\section{Minimal surfaces. Local existence, uniqueness}

\noindent A harmonic mapping is  a minimal surface, i.e. a critical point of the area functional, if and only if it is conformal. In Riemannian signature the solution of the PDE's corresponding to conformality is an additional difficulty, in Lorentzian signature it's a matter of initial data:

\begin{Theorem}
\label{min}

Assumptions as in Theorem \ref{wave} and additionally assume that the curve $k$ is such that $\langle \partial _x k_0 , \partial_x k_0 ßrangle  = - \langle k_1 , k_1 \rangle$ and $ \partial _x k_0 \perp k_1$ (note that the first requirement can be satisfied by reparametrization of the curve). Then the corresponding wave map is conformal.

\end{Theorem}

\noindent {\em Proof.} This is just a consequence of the facts that for such a curve

\begin{itemize}

\item{$||u|| = 0 = ||v||$ on the x-axis,}

\item{$||u|| = 0 = ||v||$ everywhere where the surface is defined (because of parallel transport) }

\item{therefore the pulled-back metric is a multiple of the Minkowski metric (with a possibly not overall-positive conformal factor) \qed}

\end{itemize}

\noindent Now we want to exclude sign changes of the conformal factor which is positive on the x-axis.

\begin{Theorem}
\label{gaensemarsch}
Let $M$ be time-oriented. If the curve $k: \R \rightarrow M$ satisfies $|| \partial _x k_0|| = - ||k_1|| \geq 0$, $ \partial _x k_0 \perp k_1$ and $k_1 \pm \partial_x k_0 \neq 0$ at every point, then the corresponding wave map has always

$$|| \partial_t y|| \leq 0 , \qquad || \partial_x y || \geq 0 $$ and thus is a Lorentzian minimal surface, i.e. a critical point of the area functional whose pulled-back metric is conformally equivalent to the standard Minkowski metric.

\end{Theorem}

\noindent {\em Proof.} We then have that on the x-axis always $u,v \neq 0$. As $u,v$ are parallel-transported over the surface along characteristic lines we know that both of them are nowhere equal to the zero vector and therefore both of them stay on the forward light cone because of time-orientation. So $\partial_t y = \frac{1}{2} (u+v)$ is (being a convex combination) contained in the solid forward light cone. Therefore $\p_t y$ is always causal and future-directed while $\p_x y$ is spacelike or lightlike   \qed

\bigskip

\noindent Note that however, it {\em can} happen that $u,v$ touch the light cone somewhere (and thus coincide there): Consider e.g. the case of a $k_0$ being just a circle of radius $r$ lying on the $x_1x_2$-plane in 3-dimensional Minkowski space, $k_1 = x_0$ constant, the rotational-symmetric solution $y (x,t) = (r \cos t \sin x , r \cos t \cos x , t)$ is one with cosine-shaped radius, so if the length of the string is $2 \pi$ the circle $\{ t= \pi \} $ is mapped onto a single point in $\R^{1,3}$, $\partial_x y = 0$ there.  Note that the resulting minimal surfaces will in general neither maximize nor minimize the area functional (cf. \cite{la}).

\section{Global existence}

\noindent In the following we consider only globally hyperbolic Lorentzian manifolds, i.e. $M$ is diffeomorphic to $\R \times N$ and the Lorentzian metric $g$ at a point $p=(t, x_1...x_n)$ is a direct product $-dt^2 + g^{ij} (p) dx_i \x dx_j$ in local coordinates. As additional Riemannian metric $h$ we take now the {\em flip metric to g}, i.e. $+dt^2 + g^{ij} (p) dx_i \x dx_j$. Of course, this depends on the choice of coordinates, but it does not matter for the proof.

\noindent Moreover, we shall restrict ourselves to the case of {\em closed} curves as initial data. These give rise to smooth maps from subsets of flat cylinders into spacetime because periodic initial data produce periodic solutions (as horizontal translations are conformal transformations). As one can expect, we do not get a {\em worldsheet} notion of globality. Every statement about global existence on the world-sheet in a reparametrization-invariant theory would be a pure artefact of the chosen gauge as global existence has no gauge-invariant meaning. Instead, we have a {\em target space} global existence theorem, i.e., whatever spacelike submanifold of spacetime we take, its preimage under $y$ will be a one-sphere on the cylindrical worldsheet.

\begin{Theorem}
Let $M$ be globally hyperbolic, let $k=(k_0,k_1) : \R \rightarrow TM$ be a smooth $2 \pi$-periodic curve with $k_1 \pm \partial_x k_0 \neq 0$ at every point, $|| \partial _x k_0|| = - ||k_1|| \geq 0$ and $ \partial _x k_0 \perp k_1$. Then there is a unique smooth (and $2 \pi$- periodic) map from an open subset $\Omega$ of the upper half plane $y : \{ x \geq 0 \} \rightarrow M $ which is a minimal surface and whose restriction and normal derivative at $\{ t=0 \}$ correspond to the given curve, i.e. $y(0,x) = k_0 (x)$, $\partial_t y(0,x) = k_1(x)$, and with the property that for all $p \in \partial \Omega \setminus \{ t=0 \}$

$$\lim_{x \rightarrow p} y^0 (x) = \infty ,$$

\noindent i.e. $y^{-1} ( \{ x_0 \} \times N )$ is always a smooth $2 \pi$-periodic graph over the x-axis in $\R^{1,1}$. 

\end{Theorem}

\noindent {\em Proof.} First define the maximal domain of existence of $y$ , $\Omega$, as the largest open and connected region in $\R^{1,1}$ containing the x-axis on which a solution is defined. Now assume that there is a real number $T$ such that $y^{-1}(\{T\} \times N \nsim \S^1$. As this is a closed interval, we can choose $T$ as the minimum of these zero component bounds. Now observe the curves $y^{-1} ( \{ \tau \} \times N)$ for $t < T$ on the cylinder. These are smooth closed spacelike curves on the cylinder which are therefore graphs over the x-axis such that we can write them as $(x, f_{\tau} (x) )$ with real functions $f_{\tau}$.. Now from conformality of $y$ and from the proof of Theorem \ref{gaensemarsch} one sees that $y$ always maps the backward null cone to the backward null cone in spacetime; thus $\partial_{\xi} y^0, \partial_{\eta} y ^0 > 0 $ ($y^0$ is increasing along upward characteristic lines). Thus the $f_{\tau}$ are Lipschitz functions with Lipschitz constant 1, and their limit curve corresponding to $f_T := lim_{\tau \rightarrow T} f_{\tau}$ either is situated at infinity (then $\Omega $ is the whole cylinder) or is a proper curve on the cylinder which is not necessarily smooth anymore but still Lipschitz with Lipschitz constant 1. Let us consider the second case first. The limit curve contains points of $\Omega$ and points of its complement which is closed and therefore compact, thus we can choose an $\underline{x}$ such that $f_T (\underline{x}) \notin \Omega$ and is one of the values $x$ where $f_T (x)$ is minimal among all $x$ for which $(x, f_T(x)) \notin \Omega$.

\bigskip

 Then because of the twofold minimal choice of this point the backward causal cone of $(\underline{x}, \underline{t})$ in $\R^{1,1}$ belongs entirely to the domain of existence with the exception of $(\underline{x}, \underline{t})$ itself, and in this cone we have $y^0 \leq T$. Within this cone we choose now an arbitrary characteristic triangle $\Delta_0$ with $(\underline{x}, \underline{t})$ as top point. The base side of this triangle we call $a$.
 This gives us suitable $L^1$-bounds for $||u||$ and $||v||$ in this triangle in the following way: Let $T_0 := \min \{ y^0(p) | p \in a \}$. Then we have for every point $p \in \Delta_0 \setminus (\underline{x}, \underline{t})$ the estimate

$$ \int_{c_{p, \xi}} || u || = \sqrt{2} \int_{c_{p, \xi}} u^0 \leq \sqrt{2} (T-T_0)$$

\noindent Therefore, $y (\Delta_0 ) \subset B_{\sqrt{2}(T-T_0)} (a)$. Choose the triangle $\Delta_0$ so small that this ball is contained in some $h$-ball of half the injectivity radius at its midpoint. Let $G$ be an upper bound for the operator norm of the Christoffel symbols in this ball. Now, the global upper estimate $\sqrt{2} (T-T_0 )$ for all integrals of the form $\int_{c_{\xi}} ||u||$ and  $\int_{c_{\eta}} ||v||$ together with $G$ imply global $C^0$-bounds for $||u||,||v||$ in $\Delta_0 \setminus (\underline{x}, \underline{t} ) $ because of the defining differential equation 

$$ 0 = D_{\eta} u = \partial_{\eta} u^{\gamma} (x) + \Gamma^{\gamma}_{\a \b} (y) u^{\a} (x) v^{\b} (x) .$$  

\noindent (cf. \cite{wa}). Now the length of the triangles one can add is proportional to $ \frac{1}{||u||+||v||}$, but this quantity is bounded from zero now in $\Delta_0 \setminus (\underline{x}, \underline{t})$, thus somewhere one can find a horizontal line from the left to the right side of the triangle which one can prolong still a bit over the endpoints such that it is still in the injectivity radius ball (which is possible since $\Omega$ is open) and satisfies the conditions of the local statement (as in the estimates of the local statement there was enough space left). Thus, one can add a triangle containing the point $(\underline{x}, \underline{t})$ in its interior, a contradiction.

Thus, for all $x$,  

$$ \lim_{p \rightarrow \partial \Omega \setminus \{ t= 0 \} } y^0 (p) = \infty  $$  

\bigskip

\noindent Now the only case that cannot be treated like this is when $\underline{t} = \infty$. As above, from the definition of $(\underline{x},\underline{t})$ it can be easily seen that then $y$ is defined globally on the cylinder. But then following suitable characteristic lines one can show that the whole halfplane $C^1$-converges to one point $P$ in the image for $t \rightarrow \infty$ (which again enables us to get an a priori bound on the Christoffel symbols). The proof can then be given by using Lorentzian normal coordinates at this point: Let the solution already be in the $\delta$-ball around $P$, with $\delta \leq R$. Then we write down the ODE estimates for a solution of system (\ref{vier}) compared to a solution $\tilde{u}$ of the flat equations, i.e. with all $\Gamma$ set to zero, with same initial values $\underline{u}$, along a characteristic line $c(t)$:

If $|| \Gamma^{\a}_{\b \gamma} (c(t)) \cdot v^{\gamma} (c(t)) || \leq h(t)$ then 

$$ ||u(t) - \tilde{u} (t) || \leq e^{\int_0^t h(t)}  \int_0^t e^{- \int_0^s h(s)} (   || \Gamma^{\a}_{\b \gamma} (c(s)) \cdot v^{\gamma} (c(s)) || \cdot ||\tilde{u} (s) || ) ds$$

$$ \leq e^{\int_0^t h(t)}  \int_0^t  (   || \Gamma^{\a}_{\b \gamma} (c(s)) \cdot v^{\gamma} (c(s)) || \cdot ||\tilde{u} (s) || ) ds$$

$$ \leq e^{\int_0^t G(s) ||v(s)|| ds}  \int_0^t  (  G(s) \cdot ||v(s) ||  \cdot ||\underline{u}|| ) ds$$

\noindent where $G(s)$ is an upper bound for the operator norm of the Christoffel symbols at $y(x, \tau )$ for $ \tau \geq t(c(s))$, a monotonically decreasing function with limit zero because of the choice of Lorentzian normal coordinates. Therefore

$$ ||u(t) - \tilde{u} (t) || \leq e^{G \delta}  G \underline{ u} \delta $$

\noindent Now choose a $\delta$ such that $e^{G \delta} G \delta \leq \frac{1}{2}$. Then $y^0 (c_{\xi} (l) )- y^0 (c_{\xi} (0)) \geq   l \cdot \underline{u} - l \cdot e^{G \delta} G \delta  \underline{ u} $ (the first part of the sum corresponds to the solution in flat target space while the second one is the maximal possible distance from this flat solution). Once again:
$y^0 (c_{\xi} (l) )- y^0 (c_{\xi} (0)) \geq   l \cdot \underline{u} - l \cdot e^{G \delta} G \delta  \underline{ u} \geq \frac{1}{2}l \underline{u} \geq 2 \delta$ for $l$ large enough, therefore $y^0$ would transgress $P^0 + \delta$, a contradiction  \qed  

\bigskip

\noindent Of course, one can reverse the direction of the process and solve the minimal surface equations backwards in time. Then one ends with a map $y$ from an open set of a cylinder into spacetime, and as $ y^0 (x,t)$ is a function strictly monotone in t with limits $\pm \infty$, $y$ crosses every spacelike submanifold of spacetime in the image of a one-sphere (which can, of course, degenerate to a point e.g. in the example of the previous section).

\section{Summary}

\noindent Now let us summarize the main result. The calculations above show that for any given curve with the properties described above there is exactly one minimal surface parametrized by conformal gauge, i.e. such that the pulled-back metric is conformally equivalent to the two-dimensional Minkowski metric on the cylinder. As mentioned before, every metric on the cylinder with Lorentzian ends can be brought into this gauge (cf. \cite{ku}). Therefore this initial curves are in one-two-one correspondence to {\em unparametrized} minimal surfaces. On the other hand we can give up parametrizations of the initial curves and consider them just as one-dimensional submanifolds of the total space of $CM$ where $ \pi : CM \rightarrow M$ is just the subbundle of $TM$ of all timelike tangential vectors. In summary, for $M$ being globally hyperbolic, we have established the following one-to-one correspondences:

\bigskip

\noindent closed spacelike curves in $M$ with an orthogonal timelike vector field along it, modulo reparametrizations 

\bigskip

\noindent $\leftrightarrow$ closed spacelike curves $c$ in $M$ with an orthogonal timelike vector field $V$ along it, with $\vert \vert V \vert \vert = - \vert \vert \dot{c} \vert \vert$

\bigskip

\noindent $\leftrightarrow$ minimal Lorentzian surface maps from an cylindrical open subset $O$ of a cylinder (with $\S^1 \times \{ 0 \}$ as the lower boundary of $\overline{O}$) into $M$ in conformal gauge (i.e. such that the pulled-back metric is conformally equivalent to the standard Lorentzian metric on the cylinder), modulo rigid rotations

\bigskip

\noindent $\leftrightarrow$ unparametrized minimal Lorentzian surface maps from a cylinder into $M$

\bigskip

\noindent where the arrows from below in inverse direction stand for the choice of the conformal parametrization, the restriction to the boundary, and finally forgetting the parametrization. Note that despite of the somewhat unusual notion of global solution, the notion is covariant in the sense that it does not depend on the choice of a time slice. 

\bigskip

\noindent Stability in characteristic triangles or corresponding stripes is easily proven by the same energy estimates with $w$ nonvanishing on the lower boundary this time. Thus intersecting with any spacelike Cauchy submanifold is not only a one-to-one correspondence but a diffeomorphism between Frechet manifolds.

\end{document}